\documentclass[11pt]{amsart}
\usepackage{amssymb}
\usepackage{amsmath}
\usepackage{mathrsfs}
\usepackage{amsfonts}
\usepackage{txfonts}
\usepackage{amssymb}
\usepackage{color}
\usepackage{xcolor}

\newcommand\RR{{{\mathbb R}}}

\newtheorem{theo}{Theorem}[section]
\newtheorem{lemm}[theo]{Lemma}

\newtheorem{prop}[theo]{Proposition}
\newtheorem{rema}[theo]{Remark}

\begin{document}

\title[Landau-Soft]
{Some A priori estimates for the homogeneous Landau equation with soft potentials}

\author{R. Alexandre}
\address{R. Alexandre,
\newline\indent Irenav, Arts et Metiers Paris Tech, Ecole
Navale,
\newline\indent
Lanveoc Poulmic, Brest 29290 France
\newline\indent
and
\newline\indent
Department of Mathematics, Shanghai Jiao Tong University
\newline\indent
Shanghai, 200240, P. R. China}
\email{radjesvarane.alexandre@ecole-navale.fr}

\author{J. Liao}
\address{J. Liao,
\newline\indent
School of Science, East China University of Science and Technology
\newline\indent
Shanghai, 200237, P. R. China}
\email{liaojie@ecust.edu.cn}

\author{C. Lin}
\address{C. Lin,
\newline\indent
Department of Mathematics, Shanghai Jiao Tong University
\newline\indent
Shanghai, 200240, P. R. China}
\email{chunjin.lin@gmail.com}

\subjclass[2000]{35H10, 76P05, 84C40}

\date{}

\keywords{Landau equation, soft potential, weak solutions}

\begin{abstract}
 This paper deals with the derivation of some \`a priori estimates for the homogeneous Landau equation with soft potentials. Using the coercivity of the Landau operator for soft potentials, we prove a global estimate of weak solutions in $L^2$ space without any smallness assumption on the initial data for $ -2 < \gamma <0$. For the stronger case $ -3 \leq \gamma \leq  -2$,  which covers in particular the Coulomb case, we get such a global estimate, but in some weighted $L^2$ space and under a smallness assumption on  initial data.

\end{abstract}

\maketitle

\tableofcontents

%%%%%%%%%%%%%%%%%%%%%%%%%%%%%%%%%%%%%%%%%%%%%%%%%%%%%%%%

\section{Introduction}
\setcounter{equation}{0}

The classical homogenous Landau equation (also called Fokker-Planck-Landau equation) is a common model in kinetic theory, see Chapman-Cowling \cite{CC}   and  Lifschitz-Pitaevskii \cite{LP}. This equation is obtained as a continuous approximation of the Boltzmann equation when grazing collisions prevail, see for instance \cite{AV,AB,D,G,V} for a detailed study of the limiting process, and references therein on this subject.
It describes the evolution of the (homogeneous) density function $f(t,v)$ of particles having the velocity $v\in\RR^3$ at time $t>0$:
\begin{equation}\label{i1}
\partial_t f = \partial_{v_i} \bigg\{ \int_{v_*} a_{ij} (f_* \partial_{v_j} f - f \partial_{v_{*j}} f_*  )\bigg\},
\end{equation}
 where
\begin{equation}\label{i2}
a_{ij} (z) = | z|^{\gamma +2}\Pi_{ij} (z) , \ \Pi_{ij} (z) =\delta_{ij} - {{z_iz_j}\over{|z|^2}} , \ z\neq 0.
\end{equation}

The properties of the Landau equation depend heavily on $\gamma$. It is customary to speak of hard potentials for $\gamma>0$, and soft potentials for $\gamma\in(-3,0)$. The special cases, $\gamma=0$ and $\gamma=-3$, are called the Maxwellian and Coulomb potentials, respectively. Note the fact that the more $\gamma$ is negative, the more the Landau equation is physically interesting, see Villani \cite{V2} for a detailed survey about such considerations. We refer to \cite{DV,DV1,V,V2} for more details on this equation and its physical meanings.

For a given nonnegative initial data $f_0$, we shall use the notations
\begin{equation}\label{i3}
m_0 = \int_{\RR^3}f_0(v)dv, \hspace{5mm}
e_0 = \frac{1}{2}\int_{\RR^3}f_0(v)|v|^2dv, \hspace{5mm}
H_0 = \int_{\RR^3}f_0(v) \log f_0(v)dv,
\end{equation}
for the initial mass, energy and entropy. It is classical that if $f_0\geq 0$ and $m_0,~ e_0,~ H_0$ are finite, then $f_0$ belongs to
$$
L\log L(\RR^3)= \Big\{ f\in L^1(\RR^3):  \int_{\RR^3} |f(v)|~ |\log (|f(v)|)|dv < \infty \Big\}.
$$
The solution of the Landau equation satisfies, at least formally, the conservation of mass, momentum and energy, that is, for any $t>0$,
$$
 \int_{\RR^3}f(t,v) \varphi(v) dv = \int_{\RR^3}f_0(v) \varphi(v) dv, \hspace{5mm} \varphi(v)=1,~ v,~ |v|^2/2.
$$
We also define
$$ m = \int_{\RR^3}f(t,v)  dv = m_0, \hspace{5mm}
 e = \int_{\RR^3}f(t,v) \frac{|v|^2}{2} dv = e_0.$$
Another fundamental a priori estimate is the decay of entropy, that is, the solution satisfies, at least formally, for any $t>0$,
$$
\int_{\RR^3}f(t,v) \log f(t,v) dv \leq
\int_{\RR^3}f_0(v) \log f_0(v) dv.
$$

For $s\geq 0$, we introduce classical weighted spaces as follows
$$ \| f\|_{L^1_s(\RR^3)} =   \int_{\RR^3} |f(v)|~ <v>^s dv =M_s(f),$$
$$ \| f\|^2_{L^2_s(\RR^3)} \equiv \int |f(v)|^2 <v>^{2s}  dv , $$
where $<v>:=(1+|v|^2)^{1/2}$. And we set
$$b_i = \partial_j a_{ij} (z) , \ c(z) =\partial_{ij} a_{ij}(z),$$
$$\bar a_{ij} =a_{ij} * f , \ \bar b_i = b_i * f , \ \bar c =c *f.$$
If $\gamma > -3$, we have
$$a_{ij} = \Pi_{ij} (z) | z|^{\gamma +2}, \ b_i = -2 | z|^{\gamma +2} {{z_i}\over{|z|^2}}, \ c = -2(\gamma +3) | z|^\gamma,$$
and if $\gamma =-3$, the first two formulas remain true while the third one is replaced by
$$c= -8\pi \delta_0.$$

The theory of the homogeneous Landau equation for hard potentials is studied in  great details by Desvillettes-Villani \cite{DV,DV1}, while the particular case of Maxwellian molecules $\gamma=0$ can be found in Villani \cite{V1}.

However, there are only scattered results concerning  the soft potentials. We mention the compactness properties in Lions \cite{L} and the existence of weak solutions in the inhomogeneous context by means of renormalization tools in  Villani \cite{V0} for very soft potentials, the existence of H-solution under some assumptions on initial conditions considered in Villani  \cite{V}. By using a probabilistic approach, Guerin \cite{G1} studied the existence of a measure solution for $\gamma\in(-1,0)$. Still by probabilistic approach,  Fournier-Guerin \cite{FG} studied the uniqueness and local existence of such weak solutions for soft potentials. For the Coulomb potential case $\gamma=-3$, Arsen'ev-Peskov \cite{AP} studied the local existence of weak solutions and  Fournier \cite{F} considered the local well-posedness result for such solutions. All these results give a priori estimates of solutions in some $L^p$ spaces, globally if $-2 < \gamma  <0$ and locally if $-3 \leq \gamma \leq -2$.

%\vspace{5mm}

This paper is devoted to some further a priori energy estimates by using the coercivity of the Landau operator for soft potentials, given in  Desvillettes-Villani \cite{DV}, which is stated and proved therein for $\gamma >-2$ but remains true for $\gamma \geq -3$ (at least). Our main result is the following, where here and below we use $C$ or $C_i$ to denote a generic constant.
%\vspace{2mm}

\begin{theo}\label{Th} Consider the Cauchy problem for the classical homogenous Landau equation (\ref{i1})-(\ref{i2}) with initial data $f_0\in L^1_2(\RR^3)\cap L \log L(\RR^3)$. Let the initial mass $m_0$, energy $e_0$ and entropy $H_0$ defined in (\ref{i3}) be finite. Then we have
\vspace{1mm}

\noindent {\bf 1.} Assume that $ -2 \leq \gamma <0$ and $f_0\in L^2(\RR^3)$. Then we have the following global in time a priori estimate on a weak solution in $L^2(\RR^3)$
$$
\| f(t)\|^2_{L^2(\RR^3)} \leq e^{C_2 t } \Big(\| f_0\|^2_{L^2(\RR^3)} + C_1 t \Big),
$$
where the constants $C_1$ and $C_2$ depend on  $\gamma,~ m_0,~ e_0$ and $H_0$.

\vspace{1mm}

\noindent {\bf 2.} Assume that $-3 \leq \gamma < -2$ and that $f_0\in L^2_\alpha(\RR^3)$ for some  $\alpha\geq -1-3\gamma/2$. Assume moreover that  $\|f_0\|_{L^2_\alpha(\RR^3)}$ is suitably small. Then there exists a constant $\tilde C$ depending only the entropy estimates of $f_0$ and on $\|f_0\|_{L^2_\alpha(\RR^3)}$ such that one has a global in time a priori estimate on a weak solution
$$f\in L^\infty  ([0,\infty);L^2_\alpha(\RR^3) )
\hspace{2mm} \text{ and} \hspace{2mm}
 \| f (t,\cdot)\|_{L^2_\alpha(\RR^3)}  \leq \tilde C .$$

\noindent {\bf 3.} Under assumptions stated in 1. or 2. above, one has
$$f \in L^2 (0, T ; H^1_{\alpha} (\RR^3_v )) \mbox{ for any fixed } T>0,$$
where $H_\alpha$ denotes the corresponding weighted Sobolev space.

\end{theo}

We note that from these a priori energy estimates in weighted $L^2$ spaces and similar ones for higher derivatives, eventually with different weight functions which can be obtained following the general scheme displayed below, one could get the complete existence result by using the arguments of  Desvillettes-Villani \cite{DV} and Arsen'ev-Peskov \cite{AP}. In particular, one could eventually have an immediate regularization property of solutions.

Moreover, we remark that uniqueness and convergence to equilibrium results of these weak solutions can be derived based on the works of Fournier \cite{F} and  Fournier-Guerin \cite{FG}. Note also that we decided to work in $L^2$ type spaces, but our proofs can also be adapted to more general weighted $L^p$ spaces with $1<p<+\infty$, at the expense of changing one crucial argument used in the proofs, namely Pitt's inequality, see Beckner \cite{B,B1,B2} for example. Finally, a comparison with the recent result of Fournier-Guerin \cite{FG} shows that we slightly improve their results even in the case $ \gamma >-2$ but close to $\gamma =-2$, and of course in the case $-3\leq \gamma \leq -2$, though we need a smallness assumption.

The proof of our main result above rests mainly on Pitt's inequality \cite{B,B1,B2}. However, it is possible to avoid this inequality at least in the case of not too soft potentials $\gamma \in (-2,0)$, by using standard Nash Gagliardo Nirenberg inequalities \cite{TAY} for example, and assuming enough control of moments in $L^1$, as follows from Villani \cite{V}. For example, one can show that

\begin{prop}\label{prop-radja}
Under the same hypothesis as in Theorem \ref{Th}, assume moreover that $\gamma \in (-2,0)$, that $M_\mu (t)$ is bounded by $c_\mu (1+t)$, where $\mu = {{-4\gamma (3-\gamma )}\over{3(2+\gamma )}}$. Then it follows that

$$\| f(t)\|^2_{L^2(\RR^3)} \leq C(1+t)^2 .$$
\end{prop}

Comparing with Theorem \ref{Th}, we improve on the temporal growth of this $L^2$ norm. But we do ask for many more moments: in particular, note that for $\gamma$ very close to $-2$, then we ask for almost all moments to be controlled. This point might be linked with working with $L^2$ type estimations, see last Section for further comments.

%\vspace{5mm}

The organization of the paper is as follows. Firstly, a proposition of coercivity for soft potentials is proved in Section 2, following the arguments of Desvillettes-Villani \cite{DV}, for $\gamma \geq -3$. This section is merely for the convenience of the reader since the proof follows by carefully looking to the proof in \cite{DV}.

Then in Section 3, the a priori energy estimates are carried out for the case $\gamma \in (-2,0)$ to get the global estimate of weak solutions, giving the first part of Theorem \ref{Th}.

In Section 4, we carry out the weighted energy estimates for the case $- 3< \gamma \leq -2$ to get the global estimates of weak solutions in weighted $L^2$ spaces, upon a smallness assumption on the initial data. This gives the second part of Theorem \ref{Th}, completed by Section 5. for the special case $\gamma =-3$.

In Section 6, again the same process is shown to yield local in time estimate for the case $\gamma \in(-3,-2)$, unless we can get better moment estimates in $L^1$ (that is, if the moment is uniformly bounded w.r.t time). But up to now, we have only a upper bound with a linear time growth according to Villani \cite{VT}.

Finally, Section 7 is devoted to the proof of Proposition \ref{prop-radja}.

%\vspace{5mm}

\section{Coercivity}
\setcounter{equation}{0}

This section is devoted to the proof of coercivity for soft potentials, which is an extension of hard potential case in Desvillettes-Villani \cite{DV}. In fact as mentioned to us by Desvillettes, the proof stated therein works for $\gamma >-2$ but we show that it still holds true for $\gamma \geq -3$.

\begin{prop}\label{coer} (Coercivity)
Let $\gamma\in [-3, 0)$. Let $f\in L^1_2\cap L \log L(\RR^3)$ with $m(f)=m_0$, $e( f ) \leq e_0$, $H ( f ) \leq H_0$. Then there exist a constant $C_{coer}$, explicitly computable and depending on $\gamma,~ m_0,~ e_0$ and $H_0$, such that
\begin{equation}\label{c0}
\forall \xi\in \RR^3, \hspace{8mm} \bar{a}_{ij}\xi_i\xi_j \geq C_{coer} <v>^\gamma |\xi|^2.
\end{equation}
\end{prop}

To prove the coercivity proposition, we use the same notations as in \cite{DV}, and recall the following lemma from \cite{DV}:

\begin{lemm}\label{lemma-smallset}
 Let $f \geq 0$ be a function of $L^1(\RR^3)$ such that $m ( f ) = m_0$, $e( f ) \leq e_0$, $H ( f ) \leq H_0$. Then, for all $\epsilon > 0$, there exists $\eta( \epsilon ) > 0$, depending only on $m_0, e_0, H_0$, such that for any measurable set $A \subset \RR^3$,
$$|A|\leq \eta( \epsilon ) ~ \Rightarrow ~ \int_A f \leq \epsilon,$$
where $|A|$ denotes the Lebesgue measure of $A$.
\end{lemm}

\noindent{\bf Proof of Lemma \ref{lemma-smallset}:} the arguments are taken from the nonhomogeneous case dealt with by Desvillettes \cite{D1}. But we slightly modify some of his steps, since we display an explicit expression of $\eta (\varepsilon )$ which could be required elsewhere (and which of course is not unique as regards of the proof below).

We note firstly that

$$\int f | \log f | - \int f\log f = 2\int_{f\leq 1} - f\log f$$
$$ =2 \int_{   e^{ -1 -1 { {|v|^2 }\over 2  } } \leq f \leq 1 } -f\log f +2 \int_{f \leq e^{-1- {{|v|^2}\over 2  }}} - f\log f $$

$$\leq 2m +2 e + 3(2\pi )^3\exp (1).$$

Using the decrease of entropy, and the conservation of mass and energy, it follows that

$$\int f | \log f |  \leq H_0 + 2m_0 +2 e_0 +3(2\pi )^3\exp (1) \equiv \tilde H_0.$$

Now let a fix an arbitrary set $A$. One has, for all $\delta \geq 1$

$$\int_A f = \int_{A, f\leq \delta} f + \int_{A, f\geq \delta} f$$
$$\leq \delta | A| +  (\log \delta )^{-1} \tilde H_0 .$$

Assume that $|A| \leq {{\tilde H_0} \over{\delta \log\delta}}$. Then it follows that

$$\int_A f \leq 2 \tilde H_0 {1\over{\log \delta}} .$$

We want this to be less than a fixed $\varepsilon >0$. It is enough to take the value of $\delta$ as $\delta = e^{2\tilde H_0 \varepsilon^{-1}}$. In conclusion, we have shown that setting

$$\eta (\varepsilon ) = {\varepsilon\over 2} e^{-2 \tilde H_0 \varepsilon^{-1}},$$
then it follows that

$$\int_A f \leq \varepsilon \mbox{ for any measurable set } A \mbox{ such that } | A| \leq \eta (\varepsilon ),$$
ending the proof.

\smallskip

\noindent{\bf Proof of Proposition \ref{coer}:}
Let $\xi \in \RR^3$, $| \xi | =1$, $0<\theta < {\pi\over 2}$. And set

$$D_{\theta,\xi } (v) = \Big\lbrace v_\ast \in \RR^3 :  | {{v-v_*}\over{| v-v_*|}} \cdot \xi | \geq \cos \theta \Big\rbrace,$$
which is the cone centered at $v$, of axis directed by by $\xi$ and of angle $\theta$ (see the figure in \cite{DV}).

For all $v_*\in \RR^3 \backslash D_{\theta,\xi} (v)$, we have

$$a_{ij} (v-v_*) \xi_i \xi_j = | v-v_*|^{\gamma +2} \Big(\delta_{ij} - {{(v-v_*)_i (v-v_*)_j}\over{|v-v_*|^2}} \Big) \xi_i \xi_j$$
$$ = |v-v_*|^{\gamma +2} [ 1- | {{v-v_*}\over{|v-v_*|}}\cdot\xi|^2 ] \geq |v-v_*|^{\gamma +2}\sin^2 \theta .$$
Then for all $v\in \RR^3$, $\theta \in (0, {\pi\over 2})$, $R_* >0$, we get

$$\bar a_{ij} (v) \xi_i \xi_j \geq \int_{\RR^3 \backslash D_{\theta ,\xi }(v) } dv_* f_* 1_{|v_*| \leq R_*} a_{ij} (v-v_*) \xi_i \xi_j$$

\begin{equation}\label{c1}
\geq \int_{\RR^3 \backslash D_{\theta ,\xi }(v) } dv_*  1_{|v_*|\leq R_*} |v-v_*|^{\gamma +2} f_* \sin^2 \theta.
\end{equation}

We first take care of large $|v_*|$. Let $R_*=2(e_0/m_0)^{1/2}$ and $B_*$ be the ball with center $0$ and radius $R_*$. Then
\begin{equation}\label{c2}
\int_{B_*}dv_* f_* \geq m_0(1-\frac{2e_0}{m_0R_*^2})\geq \frac{m_0}{2},
\end{equation}
and we also note that
\begin{equation}\label{c3}
| B_* \cap D_{\theta ,\xi }(v) |  \leq 2\pi R_* (|v|+R_*)^2 \tan^2\theta.
\end{equation}

We consider two cases:
\vspace{3mm}

\noindent{\bf Case 1:} $|v|\geq 2R_*.$ Note that

$$
\forall v,v_*\in \RR^3, ~~~|v|\geq 2R_*, ~~~|v_*|\leq R_*:  \hspace{5mm}
   \frac{1}{2}|v| \leq |v-v_*| \leq \frac{3}{2}|v|,
$$
from (\ref{c1}), we have
\begin{equation}\label{c4}
\bar a_{ij} (v) \xi_i \xi_j \geq\left\{
\begin{array}{ll}
(\frac{1}{2}|v|)^{\gamma+2} \sin^2 \theta \int_{B_* \backslash D_{\theta ,\xi }(v) } dv_* f_*  ~~~~~\text{for}~~~~~ \gamma\in[-2,0), \vspace{3mm} \\
(\frac{3}{2}|v|)^{\gamma+2} \sin^2 \theta \int_{B_* \backslash D_{\theta ,\xi }(v) } dv_* f_*  ~~~~~\text{for}~~~~~ \gamma\in[-3,-2).
\end{array}
\right.
\end{equation}

Now we choose $\theta>0$ such that
$$ \tan^2 \theta = \min\{\frac{2\eta(\frac{m_0}{4})}{9\pi R_*|v|^2}, ~1 \},$$
so according to (\ref{c3}), we have
$$
| B_* \cap D_{\theta ,\xi }(v) |  \leq 2\pi R_* (\frac{3}{2}|v|)^2 \tan^2\theta \leq \eta(\frac{m_0}{4})
$$
thus
$$
\int_{B_* \cap D_{\theta ,\xi }(v)} dv_* f_*  \leq \frac{m_0}{4},
$$
and then from (\ref{c2}) and (\ref{c4}) we have
$$
\bar a_{ij} (v) \xi_i \xi_j \geq\left\{
\begin{array}{ll}
(\frac{1}{2}|v|)^{\gamma+2} \cos^2 \theta ~ \min\{\frac{2\eta(\frac{m_0}{4})}{9\pi R_*|v|^2}, 1 \} (\frac{m_0}{2}-\frac{m_0}{4}) \geq  \frac{m_0}{32}  (\frac{1}{2}|v|)^{\gamma} \min\{\frac{2\eta(\frac{m_0}{4})}{9\pi R_*}, 4R_*^2 \}   ~~~~~\text{for}~~~~~ \gamma\in[-2,0), \vspace{3mm} \\
(\frac{3}{2}|v|)^{\gamma+2} \cos^2 \theta ~ min\{\frac{2\eta(\frac{m_0}{4})}{9\pi R_*|v|^2}, 1 \} (\frac{m_0}{2}-\frac{m_0}{4}) \geq  \frac{9m_0}{32}  (\frac{3}{2}|v|)^{\gamma} \min\{\frac{2\eta(\frac{m_0}{4})}{9\pi R_*},  4R_*^2 \}  ~~~~~\text{for}~~~~~ \gamma\in [-3,-2),
\end{array}
\right.
$$
which is also
\begin{equation}\label{c5}
\bar a_{ij} (v) \xi_i \xi_j \geq c |v|^{\gamma} \geq c<v>^{\gamma} \hspace{5mm}  \text{for}  \hspace{5mm}    \gamma\in[-3,0),~ |v| \geq 2R_*.
\end{equation}
where $c$ is a constant depending on $\gamma,~ m_0,~  e_0$ and $H_0$.

\vspace{3mm}
\noindent{\bf Case 2:} $|v|\leq 2R_*.$ Note that when $|v_*|\leq R_*$, we have $| v-v_* | \leq 3R_*$ thus
$$
\int_{B_* \backslash D_{\theta ,\xi }(v) } dv_* f_* |v-v_*|^{\gamma +2} \geq \int_{B_* \backslash D_{\theta ,\xi }(v) } dv_* f_* |v-v_*|^{\gamma +3}/(3R_*)
$$

$$
\geq \frac{[\frac{3}{4\pi} \eta(\frac{m_0}{8}) ]^{ \frac{\gamma+3}{3}}}{3R_*}
\int_{B_* \backslash D_{\theta ,\xi }(v) } dv_* f_*
1_{ |v-v_*|\geq [\frac{3}{4\pi} \eta(\frac{m_0}{8})  ]^{1/3} }.
$$
Note that
$$
\int_{B_* \backslash D_{\theta ,\xi }(v) } dv_* f_*
1_{ |v-v_*|\geq [\frac{3}{4\pi} \eta(\frac{m_0}{8})  ]^{1/3} }
= \int_{B_* \backslash D_{\theta ,\xi }(v) } dv_* f_*
 - \int_{B_* \backslash D_{\theta ,\xi }(v) } dv_* f_*
1_{ |v-v_*|\leq [\frac{3}{4\pi} \eta(\frac{m_0}{8})  ]^{1/3} }
$$
$$
\geq \int_{B_* } dv_* f_*
-\int_{B_* \cap D_{\theta ,\xi }(v) } dv_* f_* - \frac{m_0}{8},
$$
and we know the first term is greater than $m_0/2$ from (\ref{c2}) . For the second term to be less than $m_0/8$, we expect
$$
| B_* \cap D_{\theta ,\xi }(v) |  \leq 2\pi R_* (3R_*)^2 \tan^2\theta \leq \eta(\frac{m_0}{8}),
$$
which requires
$$ \tan^2 \theta = \min\{\frac{\eta(\frac{m_0}{8})}{18\pi R_*^3}, ~1 \},$$
thus from (\ref{c1}) we have
\begin{equation}\label{c6}
\bar a_{ij} (v) \xi_i \xi_j \geq \frac{[\frac{3}{4\pi} \eta(\frac{m_0}{8}) ]^{ \frac{\gamma+3}{3}}}{3R_*}  \frac{m_0}{4} \sin^2 \theta \geq \frac{[\frac{3}{4\pi} \eta(\frac{m_0}{8}) ]^{ \frac{\gamma+3}{3}}m_0}{24R_*} \min\{\frac{\eta(\frac{m_0}{8})}{18\pi R_*^3}, ~1 \},~~~~~\text{for}~~~~~  |v| \leq 2R_*.
\end{equation}
Estimates (\ref{c5}) and (\ref{c6}) together ensure the validity of (\ref{c0}).   $\hfill\square$

\begin{rema}
From the proof, we can see that actually the coercivity proposition holds for all $\gamma<0$.
\end{rema}

\vspace{2mm}

\section{The case $\gamma \geq -2$: energy estimates}
\setcounter{equation}{0}

We multiply the equation (\ref{i1}) by $f$ and integrate to get
\begin{equation}\label{e1}
{d\over{dt}} {1\over 2} \| f\|^2_{L^2} +\int\int a_{ij} (v-v_*) f_* \partial_{v_i} f\partial_{v_j} f =
{1\over 2}\int\int a_{ij} (v-v_*) \partial_{v_{*j}} f_* \partial_{v_i} f^2.
\end{equation}

The second term on the l.h.s. can be bounded below by using the coercivity property (\ref{c0}) thus
\begin{equation}\label{e2}
\int\int a_{ij} (v-v_*) f_* \partial_{v_i} f\partial_{v_j} f
\geq C_{coer} \int <v>^{\gamma} |\nabla_v f|^2 dv.
\end{equation}

For the nonlinear term arising on the on the r.h.s., we have
\begin{equation}\label{e3}
{1\over 2}\int\int a_{ij} (v-v_*) \partial_{v_{*j}} f_* \partial_{v_i} f^2 = (\gamma +3) \int_{v_*} \int_v |v-v_*|^\gamma f_* f^2   \equiv A_1+A_2,
\end{equation}
where
$$A_1 = (\gamma +3) \int_{v_*} \int_{|v-v_*| \geq R} |v-v_*|^\gamma f_* f^2$$
and
$$A_2 = (\gamma +3) \int_{v_*} \int_{|v-v_*| \leq R} |v-v_*|^\gamma f_*. f^2$$

For the first term in (\ref{e3}), since $\gamma <0$, we have
\begin{equation}\label{e4}
A_1 \leq (\gamma +3)R^\gamma  \int_{v_*} \int_{|v-v_*| \geq R} f_* f^2 \leq (\gamma +3)R^\gamma m \| f\|^2_{L^2}.
\end{equation}

The second term can be estimated as follows
$$A_2 = (\gamma +3) \int_{v_*} \int_{|v-v_*| \leq R} |v-v_*|^\gamma f_* f^2$$
$$ =  (\gamma +3) \int_{v_*} \int_{|v-v_*| \leq R} |v-v_*|^\gamma <v_*>^{-\gamma}f_* (<v>^{\gamma /2} f)^2  <v_*>^\gamma <v>^{-\gamma}.$$
Since $   <v>^{-\gamma} / <v_*>^{-\gamma}   \leq c_\gamma (1+R^2)^{-\gamma /2}$, we get
$$A_2 \leq c_\gamma (\gamma +3)(1+R^2)^{-\gamma /2}  \int _{v_*} <v_*>^{-\gamma}f_*  \int _v |v-v_*|^\gamma (<v>^{\gamma /2} f)^2 .$$
We use Pitt's inequality \cite{B,B1,B2} to get that
$$\int _v |v -v_*|^\gamma (<v>^{\gamma /2} f)^2  \leq c_{pitt} \int_\xi | \xi|^{-\gamma }| \widehat{<v>^{\gamma/2}f} (\xi )|^2.$$
In order to use Pitt's inequality, we need that $\gamma \in (-3, 0)$.

Recalling that
$$M_{-\gamma} (t) = \int _{v_*} <v_*>^{-\gamma}f_* ,$$
we get
$$A_2 \leq 3c_\gamma (\gamma +3)(1+R^2)^{-\gamma /2}  M_{-\gamma} (t) c_{pitt} \int_\xi | \xi|^{-\gamma }| \widehat{<v>^{\gamma /2} f} (\xi )|^2.$$
For any $R_*$, write
$$\int_\xi | \xi|^{-\gamma }| \widehat{<v>^{\gamma /2} f} (\xi )|^2 \leq \int_{| \xi | \leq R_*} | \xi|^{-\gamma }| \widehat{<v>^{\gamma /2} f} (\xi )|^2 + \int_{|\xi | \geq R_*} | \xi|^{-\gamma }| \widehat{<v>^{\gamma /2} f} (\xi )|^2,$$
and using the fact that $f$ is in $L^1$, we get
 $$\int_\xi | \xi|^{-\gamma }| \widehat{<v>^{\gamma /2} f} (\xi )|^2 \leq m^2\int_{| \xi | \leq R_*} | \xi|^{-\gamma}  + \int_{|\xi | \geq R_*} | \xi|^{-\gamma }| \widehat{<v>^{\gamma /2} f} (\xi )|^2$$
Now we need to assume $-\gamma  \leq 2$:
$$\int_\xi | \xi|^{-\gamma }| \widehat{<v>^{\gamma /2} f} (\xi )|^2 \leq R_*^{-\gamma  +3} m^2 + {1\over{R_*^{\gamma  +2}}}\int_{|\xi | \geq R_*} | \xi|^{-\gamma }| \xi|^{\gamma +2}|\widehat{<v>^{\gamma /2} f} (\xi )|^2$$
$$\leq R_*^{-\gamma  +3} m^2 + c_{parseval} {1\over{R_*^{\gamma  +2}}}\ \| \nabla_v [   <v>^{\gamma /2} f] \|^2_{L^2}$$
$$\leq \max (m^2, c_{parseval} ) [ R_*^{-\gamma  +3} + {1\over{R_*^{\gamma  +2}}}\ \| \nabla_v [   <v>^{\gamma /2} f] \|^2_{L^2} ]$$
Optimizing w.r.t. $R_*$, we find
$$\int_\xi | \xi|^{-\gamma }| \widehat{<v>^{\gamma /2} f} (\xi )|^2 \leq 2\max (m^2, c_{parseval}) \| \nabla_v [   <v>^{\gamma /2} f] \|^{{-\gamma +3}\over{5}}_{L^2}.$$

All in all, we have obtained, for   $\gamma \geq -2$,
$$A_2 \leq 3c_\gamma (\gamma +3)(1+R^2)^{-\gamma /2}   c_{pitt} 2\max (m^2, c_{parseval}) M_{-\gamma} (t) \| \nabla_v [   <v>^{\gamma /2} f] \|^{{-\gamma +3}\over{5}}_{L^2}.$$
Fix a small $\varepsilon >0$. Then
$$A_2 \leq \varepsilon^{{\gamma -6}\over{5}}3c_\gamma (\gamma +3)(1+R^2)^{-\gamma /2}   c_{pitt} 2\max (m^2, c_{parseval}) M_{-\gamma} (t) \varepsilon^{{-\gamma +6}\over{5}}\| \nabla_v [   <v>^{\gamma /2} f] \|^{{-\gamma +6}\over{5}}_{L^2}.$$
And apply Young's inequality for product with $p = {{10}\over{-\gamma +6}}$,  $q= {{10}\over{4+\gamma}}$, we obtain
$$A_2 \leq {1\over q} \bigg\{ \varepsilon^{{\gamma -6}\over{5}}3c_\gamma (\gamma +3)(1+R^2)^{-\gamma /2}   c_{pitt} 2\max (m^2, c_{parseval}) M_{-\gamma} (t) \bigg\}^q $$
$$+  {1\over p} \bigg\{ \varepsilon^{{-\gamma +6}\over{5}}\| \nabla_v [   <v>^{\gamma /2} f] \|^{{-\gamma +6}\over{5}}_{L^2} \bigg\}^p,$$
which is also
\begin{equation}\label{e5}
A_2 \leq {{4+\gamma}\over{10}} \bigg\{ \varepsilon^{{\gamma -6}\over{5}}3c_\gamma (\gamma +3)(1+R^2)^{-\gamma /2}   c_{pitt} 2\max (m^2, c_{parseval}) M_{-\gamma} (t) \bigg\}^{{4+\gamma}\over{10}}
\end{equation}
$$ +  {{-\gamma +6}\over{10}} \varepsilon^2 \| \nabla_v [  <v>^{\gamma /2} f] \|^2_{L^2} .$$

Note that there are two terms in $A_2$ above: $M_{-\gamma}$ and  $\| \nabla_v [  <v>^{\gamma /2} f] \|^2_{L^2} $. First, for $\gamma \in [-2,0)$,
\begin{equation}\label{e6}
M_{-\gamma} = \int <v>^{-\gamma}f  \leq  \int (1+v^2)f =m+e.
\end{equation}

Second, we consider $\| \nabla_v [  <v>^{\gamma /2} f] \|^2_{L^2} $. Note that
$$  \nabla_v [  <v>^{\gamma /2} f] = \dfrac{\gamma}{2}  <v>^{\gamma /2-2} fv + <v>^{\gamma /2} \nabla_vf.$$
Now, since $\gamma <0$, we have
$$ \| \nabla_v [  <v>^{\gamma /2} f] \|^2_{L^2} \leq 2\bigg( \|\dfrac{\gamma}{2}  <v>^{\gamma /2-2} fv \|^2_{L^2} + \|<v>^{\gamma /2} \nabla_vf \|^2_{L^2} \bigg)$$
\begin{equation}\label{e7}
\leq   \frac{\gamma^2}{2} \|f \|^2_{L^2} + 2  \|<v>^{\gamma /2} \nabla_vf \|^2_{L^2}
\end{equation}

Recall the constant $C_{coer}$ which appears in the coercive inequality (\ref{e2}). Then we can choose $\varepsilon$ such that
$${{-\gamma +6}\over{5}} \varepsilon^2 = {1\over 2} C_{coer},$$
and then combining all the above results, we get
\begin{equation}\label{e8}
{d\over{dt}} \| f\|^2_{L^2}  +  C_{coer} \|<v>^{\gamma/2} \nabla_v f \|^2_{L^2}  \leq C_1 + C_2 \|f \|^2_{L^2},
\end{equation}
%or, by slightly change the proof, we can get
%\begin{equation}\label{e9}
%{d\over{dt}} \| f\|^2_{L^2}  +  C_{coer} \|<v>^{\gamma/2}  f \|^2_{\dot H_1}  \leq C_1 + C_2 \|f \|^2_{L^2}.
%\end{equation}
or, we just simply have
\begin{equation}\label{e10}
{d\over{dt}} \| f\|^2_{L^2}    \leq C_1 + C_2 \|f \|^2_{L^2},
\end{equation}
and therefore, by directly using Gronwall's inequality, we get
\begin{equation}\label{e11}
\| f(t)\|^2_{L^2} \leq e^{C_2 t } \Big(\| f_0\|^2_{L^2} + C_1 t \Big),
\end{equation}
then we have the first part of Theorem \ref{Th}.

\begin{rema}

1. By repeating the same process for higher derivatives, one can get the global existence of weak solutions. However, the bound depends on time, but the result does not require any assumption of smallness on the initial data. This is compatible with the works of Fournier-Guerin \cite{FG} where they have global existence in that case too, though in different $L^p$ spaces.

2. We work with the usual $L^2$ space but one can easily adapt our arguments for weighted $L^2$ spaces. This is done for example in the next section when $-3<\gamma <-2$. The same remark also applies for estimation of higher derivatives as well.
\end{rema}

%\vspace{10mm}

\vspace{2mm}

%\newpage

\section{The case $-3 <\gamma <-2$: weighted energy estimates }
\setcounter{equation}{0}

We carry out the weighted energy estimates in this section for $-3 <\gamma < -2$. Of course one can consider the general case of $\gamma \in (-3,0)$ but recall that we have already good estimates from the previous section.

We want to estimate $g= <v>^\alpha f$ in $L^2$, and we assume that $\alpha \geq -1-3/2 \gamma$, see below for the final arguments, explaining this value of the weight.

Multiplying the Landau equation by $<v>^{\alpha}$ and setting $g=<v>^{\alpha} f$, we have

$$\partial_v g =  \partial_{v_i} \bigg\{ \int_{v_*} a_{ij} (f_* <v>^{\alpha} \partial_{v_j} f - g\partial_{v_{*j}} f_*  )\bigg\}$$
$$- \bigg\{ \int_{v_*} a_{ij} (f_* \partial_{v_j} f - f \partial_{v_{*j}} f_*  )\bigg\}  \partial_{v_i} <v>^{\alpha}$$
$$ =  \partial_{v_i} \bigg\{ \int_{v_*} a_{ij} (f_* \partial_{v_j} g - g\partial_{v_{*j}} f_*  )\bigg\}$$
$$ - \partial_{v_i}  \bigg\{ \int_{v_*} a_{ij} f_* f\partial_{v_j}  <v>^{\alpha}  \bigg\}$$
$$- \bigg\{ \int_{v_*} a_{ij} (f_* \partial_{v_j} f - f \partial_{v_{*j}} f_*  )\bigg\}  \partial_{v_i} <v>^{\alpha}.$$
Multiplying by $g$ and integrating, we get
$${d\over{dt}} \| g\|^2_{L^2} =
  -  \int_v \int_{v_*} a_{ij} f_* \partial_{v_j}  g \partial_{v_i} g  + \int_v \int_{v_*} a_{ij}  \partial_{v_{*j}} f_*  g \partial_{v_i} g $$
\begin{equation}\label{w1}
+ \int_v\int_{v_*} a_{ij} f_* f   \partial_{v_i} g ~\partial_{v_j}  <v>^{\alpha}
\end{equation}
$$- \int_v \int_{v_*} a_{ij} f_* g \partial_{v_j} f \partial_{v_i} <v>^{\alpha}  + \int_v \int_{v_*} a_{ij}  fg \partial_{v_{*j}} f_*  )  \partial_{v_i} <v>^{\alpha}$$

$$= -I + II    + III  -IV +V.$$

We then estimate each of these terms. Here $I$ can be controlled through the coercivity estimation (\ref{c0}). More precisely, we have

\begin{equation}\label{w2}
I \geq C_{coer} \int_ v <v>^\gamma | \nabla_v g|^2 dv.
\end{equation}
By a similar argument as for \eqref{e7} with $f$ replaced by $g$, we have
$$
 \int_ v <v>^\gamma | \nabla_v g|^2 dv \geq
 \frac{1}{2}   \int_v  | \nabla_v [ <v>^{\gamma /2} g ] |^2 dv - \frac{\gamma^2}{4} \int_v  g^2 dv  ,
$$
and thus (\ref{w2}) can be rewritten as
\begin{equation}\label{w3}
I \geq \frac{C_{coer}}{2}  \int_v  | \nabla_v [ <v>^{\gamma /2} g ] |^2 -CC_{coer} \int_v  g^2.
\end{equation}

Now, as in the previous section, the term $II$ is
\begin{equation}\label{w13}
II = (\gamma +3)  \int \int |v-v_*|^\gamma f_* g^2 \equiv (\gamma +3) [ A+B ],
\end{equation}
where
\begin{equation}\label{w14}
A=  \int \int  1_{|v-v_*| \geq \varepsilon} |v-v_*|^\gamma f_* g^2
\leq m C \varepsilon^\gamma \| g\|^2_{L^2}
\end{equation}
and
\begin{equation}\label{w15}
B = \int \int 1_{|v-v_* |\leq \varepsilon}|v-v_*|^\gamma f_* g^2.
\end{equation}
We further decompose $B$ over the sets $\{f\leq f_*\}$ and $\{f\geq f_*\}$ to get
 $$B \leq \int \int 1_{|v-v_* |\leq \varepsilon}|v-v_*|^\gamma <v>^{2\alpha} f^3 + \int \int 1_{|v-v_* |\leq \varepsilon}|v-v_*|^\gamma <v>^{2\alpha} f_*^3,$$
and then after direct computation we obtain
\begin{equation}\label{w16}
B\leq C \varepsilon^{3+\gamma} (1+\varepsilon^2)^\alpha \int <v>^{2\alpha /3} f^3.
\end{equation}

We need to control $<v>^{2\alpha /3} f$ in $L^3$ by using a control of $<v>^\alpha f$ and of $\nabla [ <v>^{\gamma /2 +\alpha } f]$ in $L^2$, and we will do it by applying Holder's inequality \cite{BRE}.

We write $p_1=p_2=p_3=9$ and
$$<v>^{2\alpha /3}f =   f^{1\over p_1}  [<v>^{\alpha } f ]^{2\over p_2}  [ <v>^{\gamma /2 +\alpha } f]^{6\over p_3}
\cdot <v>^{2\alpha /3 - 2/9 \alpha - 2/3 (\gamma /2 +\alpha )} .$$
The last exponent is $-2/9 \alpha -\gamma /3$. Raised to the power $9$ this is $-2\alpha - 3\gamma$. We ask this number to be less than $2$: $\alpha \geq -1 - 3/2 \gamma$.

In conclusion we write, with $\alpha \geq -1 - 3/2 \gamma$,

$$<v>^{2\alpha /3}f \leq   [f<v>^2]^{1\over p_1}  [<v>^{\alpha } f ]^{2\over p_2}  [ <v>^{\gamma /2 +\alpha } f]^{6\over p_3}. $$

Finally we obtain the interpolation inequality
\begin{equation}\label{w17}
\int <v>^{2\alpha} f^3   \leq C  [ \int <v>^{2} f ]^{1\over 3}  \| <v>^\alpha f\|^{2\over 3}_{L^2}  \int | \nabla (<v>^{\alpha + \gamma /2} f )|^2 .
\end{equation}

In conclusion, we combine (\ref{w13})-(\ref{w17}) to get
\begin{equation}\label{w4}
II \leq m C \varepsilon^\gamma \| g\|^2_{L^2} + C \varepsilon^{3+\gamma} (1+\varepsilon^2)^\alpha   [ \int <v>^2 f ]^{1\over 3}  \| <v>^\alpha f\|^{2\over 3}_{L^2}  \int | \nabla (<v>^{\alpha + \gamma /2} f )|^2 .
\end{equation}

We now analyze the term $III$ which is given by

$$III = \int_v\int_{v_*} a_{ij} f_* f~   \partial_{v_i} g ~\partial_{v_j}  <v>^{\alpha} .$$

We have immediately that
{
$$  III  = \int_v\int_{v_*} a_{ij} f_* f^2~ \partial_{v_i}  <v>^{\alpha} \partial_{v_j}  <v>^{\alpha} + \int_v\int_{v_*} a_{ij} f_* f \partial_{v_i} f <v>^{\alpha} \partial_{v_j}  <v>^{\alpha}
$$
$$=   \int_v\int_{v_*} a_{ij} f_* f^2~ \partial_{v_i}  <v>^{\alpha} \partial_{v_j}  <v>^{\alpha} +
\frac{1}{2} \int_v\int_{v_*} a_{ij} f_*  \partial_{v_i} f^2 <v>^{\alpha} \partial_{v_j}  <v>^{\alpha}
$$
$$=  \frac{1}{2} \int_v\int_{v_*} a_{ij} f_* f^2~ \partial_{v_i}  <v>^{\alpha} \partial_{v_j}  <v>^{\alpha} -
\frac{1}{2} \int_v\int_{v_*} \partial_{v_i}a_{ij} f_* f^2 <v>^{\alpha} \partial_{v_j}  <v>^{\alpha}
$$
$$ -\frac{1}{2} \int_v\int_{v_*}a_{ij} f_* f^2 <v>^{\alpha} \partial_{v_i} \partial_{v_j}  <v>^{\alpha},$$
}
then
$$
| III|  \leq C \int_v\int_{v_*} |v-v_*|^\gamma f_* f^2 <v>^{2\alpha -2},
$$
(see the next section for similar arguments).

Since $f <v>^{\alpha } =g$, this term can be controlled like $B$ in (\ref{w15}), thus controlled by $II$,  and we can absorb $III$ and $II$ together to get that
$$II +III \leq m C \varepsilon^\gamma \| g\|^2_{L^2} + C \varepsilon^{3+\gamma} (1+\varepsilon^2)^\alpha [ \int <v>^{2} f ]^{1\over 3}  \| <v>^\alpha f\|^{2\over 3}_{L^2}  \int | \nabla (<v>^{\alpha + \gamma /2} f )|^2 .$$

Next for $IV$, we have
$$IV = \int_v \int_{v_*} a_{ij} f_* g \partial_{v_j} f \partial_{v_i} <v>^{\alpha} ={1\over 4}  \int_v \int_{v_*} a_{ij} f_*  \partial_{v_j} f^2 \partial_{v_i} <v>^{2\alpha} $$
$$= - {1\over 4}  \int_v \int_{v_*} \partial_{v_j}  a_{ij} f_*  f^2 \partial_{v_i} <v>^{2\alpha} - {1\over 4}  \int_v \int_{v_*} a_{ij} f_* f^2  \partial_{v_j}  \partial_{v_i} <v>^{2\alpha}$$
$$\leq C \int\int |v-v_*|^{\gamma +1} f_*f^2 <v>^{2\alpha -1} + C\int \int |v-v_*|^\gamma f_*f^2 <v>^{2\alpha -2},$$
and so again we can absorb it with earlier terms. One can see that $V$ is also similar so all in all

$$II +|III| +| IIV|+ |V|  \leq m C \varepsilon^\gamma \| g\|^2_{L^2} $$
$$+ C \varepsilon^{3+\gamma} (1+\varepsilon^2)^\alpha [ \int <v>^{2} f ]^{1\over 3}  \| <v>^\alpha f\|^{2\over 3}_{L^2}  \int | \nabla (<v>^{\alpha + \gamma /2} f )|^2 .$$

By combining the above estimations, the final conclusion is that

\begin{equation}\label{w5}
 {d\over{dt}} \| g\|^2_{L^2} +  C_{coer}  \int_v  | \nabla_v [ <v>^{\gamma /2} g ] |^2 \leq
\end{equation}
$$CC_{coer} \| g\|^2_{L^2}   +m C \varepsilon^\gamma \| g\|^2_{L^2} + C \varepsilon^{3+\gamma} (1+\varepsilon^2)^\alpha [ \int <v>^{2} f ]^{1\over 3}  \| g\|^{2\over 3}_{L^2}  \int | \nabla (<v>^{\gamma /2} g )|^2 .$$

Set $E= \int <v>^2 f dv = m+e$ (the summation of mass and energy) which is bounded uniformly in time. Setting
$$X = \| g\|^2_{L^2},   \hspace{5mm}  \| <v>^{\gamma /2} g \|^2_{\dot H_1} = \int | \nabla (<v>^{\gamma /2} g )|^2, $$
the above inequality (\ref{w5}) reads as

$$ {d\over{dt}} X+ C_{coer} \| <v>^{\gamma /2} g \|^2_{\dot H_1} \leq $$
$$CC_{coer} X  +m C \varepsilon^\gamma X + C \varepsilon^{3+\gamma} (1+\varepsilon^2)^\alpha E^{1\over 3}  X^{1\over 3} \| <v>^{\gamma /2} g \|^2_{\dot H_1} ,$$
which can be also written under the form
 \begin{equation}\label{w6}
 {d\over{dt}} X \leq - \| <v>^{\gamma /2} g \|^2_{\dot H_1}  \bigg\{  C_{coer} - C \varepsilon^{3+\gamma} (1+\varepsilon^2)^\alpha E^{1\over 3}  X^{1\over 3} \bigg\} +[CC_{coer}   +m C \varepsilon^\gamma ] X .
\end{equation}

We want to proceed as in Toscani \cite{T}. However, we have a major trouble in that the moment of order $s$ of $f$ in $L^1$ are not known to be uniformly bounded w.r.t. time, see Villani \cite{VT}. This means that using Nash's inequality as in Toscani \cite{T} at that point involves a lower bound which decays in time, and so can be very small for large time. Since we want to get global solutions, we are going to use Pitt's inequality instead of Nash's inequality.

Assume that at time $t$, we have for some $\delta \geq 0$
$$ C_{coer} - C \varepsilon^{3+\gamma} (1+\varepsilon^2)^\alpha E^{1\over 3}  X^{1\over 3}  \geq \delta, $$
that is
\begin{equation}\label{cond}
C \varepsilon^{3+\gamma} (1+\varepsilon^2)^\alpha E^{1\over 3}  X^{1\over 3}  \leq  C_{coer} -\delta,
\end{equation}
which is to be used for $\delta < C_{coer}$.

Then from (\ref{w6}) we obtain
\begin{equation}\label{w7}
{d\over{dt}} X \leq - \| <v>^{\gamma /2} g \|^2_{\dot H_1} \delta +[CC_{coer}   +m C \varepsilon^\gamma ] X .
\end{equation}

Pitt's inequality tells us that
\begin{equation}\label{w8}
\| <v>^{\gamma /2} g \|^2_{\dot H_1}  \geq C_{pitt} \int {{<v>^\gamma g^2}\over{|v|^2}} dv.
\end{equation}
We are going to show a lower bound for the r.h.s. of this inequality.

We start from

$$\int f (v)dv \leq  \int_{|v|\leq R} f(v) dv + {1\over R^2} e ,$$
which is also

$$\int f (v)dv \leq  \int_{|v|\leq R} <v>^{-\alpha}g(v) dv + {1\over R^2} e $$
$$= \int_{|v|\leq R} <v>^{-\alpha}<v>^{-\gamma /2} |v| <v>^{\gamma /2} |v|^{-1} g(v) dv + {1\over R^2} e ,$$
by using Cauchy-Schwartz inequality, we obtain
\begin{equation}\label{w9}
\int f (v)dv \leq   [ \int_{|v|\leq R} <v>^{-2\alpha}<v>^{-\gamma } |v|^2 ]^{1\over 2}   [  \int_{|v|\leq R}   <v>^{\gamma } |v|^{-2} g^2(v) dv   ]^{1\over 2} + {1\over R^2} e .
\end{equation}
Recall that we have $\alpha \geq -1-3/2 \gamma$, so the exponent
% {\red (We donot need $\gamma\leq -2$ here!!!)}
$$\tilde{\alpha}:= -2\alpha -\gamma  + 2 \leq   2 \gamma +4  \in  (-2,4) \hspace{5mm} \text{when} \hspace{5mm} \gamma\in (-3,0).$$
Thus we can estimate the upper bound inside the first integral in (\ref{w9}) to get
$$\int f (v)dv \leq   CR^{7/2}  [  \int   <v>^{\gamma} |v|^{-2} g^2(v) dv   ]^{1\over 2} + {1\over R^2} e .$$
Omitting the constant $C$, this is of the form
$$\int f (v)dv \leq   R^{7/2}  A+ {1\over R^2} e
 \hspace{2mm} \text{with}  \hspace{2mm}
  A= [  \int   <v>^{\gamma} |v|^{-2} g^2(v) dv   ]^{1\over 2}, $$
and we choose $R$ such that
$$R^{7/2}  A = {1\over R^2} e ,$$
that is
$$  R  = e^{2/11} A^{-2/11} .$$
It follows that
$$m= \int f (v)dv \leq   2  R^{7/2}A = 2 e^{7/11} A^{4/11},
 $$
and furthermore we get
$$
A \geq C m^{11/4} e^{-7/4},
$$
which is also
\begin{equation}\label{w10}
\int  <v>^{\gamma } |v|^{-2} g^2(v) dv  \geq  C m^{11/2}  e^{-7/2}.
\end{equation}

We can now go back to our differential inequality (\ref{w7}): by Pitt's inequality (\ref{w8}) and (\ref{w10}), we have
\begin{equation}\label{w11}
{d\over{dt}} X \leq - C_{pitt} C m^{11/2}  e^{-7/2} \delta +[CC_{coer}   +m C \varepsilon^\gamma ] X =F(X).
\end{equation}

Recall that we have assumed (\ref{cond}), that is
$$\delta <  C_{coer}$$
and that we want
$$C \varepsilon^{3+\gamma} (1+\varepsilon^2)^\alpha E^{1\over 3}  X^{1\over 3}  \leq  C_{coer} -\delta ,$$
that is
$$
  X\leq  [ C_{coer} -\delta  ]^3 C^{-3} \varepsilon^{- 3(3+\gamma )} (1+\varepsilon^2)^{-3\alpha } / E \equiv \tilde X .
$$
Now, let $X_{eq}$ be the zero of the function $F$ defined in (\ref{w11}).
Assume that
\begin{equation}\label{w12}X(0) \leq \bar X\equiv \min \lbrace \tilde X, X_{eq} \rbrace.\end{equation}
Then, in view of the form of the differential inequality (\ref{w11}) and the behavior of the function $F$, it follows that for all $t>0$
 $$X(t) \leq \bar X.$$
Thus we have obtain a global bound for the weighted $L^2$ norm of $f$, uniformly in time, that is, we get the second part of Theorem \ref{Th}.

\begin{rema} Note that $\tilde X$ can be as large as we want, since this quantity depends on negative powers of $\varepsilon$, which is a free parameter that we can take small.  However, taking such a small $\varepsilon$, wee that $X_{eq}$ which is given by

$$X_{eq} = {{ C_{pitt} C m^{11/2}  e^{-7/2} \delta }\over{[CC_{coer}   +m C \varepsilon^\gamma ] }}$$
is going to be small. Therefore by choosing $\varepsilon$ sufficiently small, we can therefore assume that

$$\bar X =X_{eq}$$
and it follows that we will have $X(t) \leq X_{eq}$. 

At this point, it is important to recall that for any function $f$, again using the same notation as above, one has the following interpolation inequality

$$X \geq C m^{7/2} e^{-3/2}$$

Thus we should have

$$ C m^{7/2} e^{-3/2} \leq {{ C_{pitt} C m^{11/2}  e^{-7/2} \delta }\over{[CC_{coer}   +m C \varepsilon^\gamma ] }}$$
that is

$$ {{ C_{pitt} C m^{2}  e^{-2} \delta }\over{[CC_{coer}   +m C \varepsilon^\gamma ] }} \geq \tilde C .$$

Now we note, in view of previous results on coercivity that an upper bound for $C_{coer}$ is given by 
$$\max \lbrace C_1mR_*^2 , C_2  {{( Cm )^{(\gamma +3)/3}}\over R_* }\rbrace$$

We choose the value of $R_*$ such that these two terms are equal, getting an upper bound like $m^{ (9-2\gamma )/9}$. 

Then it is enough to ask for

$$ {{ C_{pitt} C m^{2}  e^{-2} \delta }\over{[C m^{ (9-2\gamma )/9}  +m C \varepsilon^\gamma ] }} \geq \tilde C .$$

Then we choose a smaller $\varepsilon$ so that the second term on the denominator is bigger than the first one, so we are led to ask for

$$ {{ C_{pitt} C m^{2}  e^{-2} \delta }\over{[2m C \varepsilon^\gamma ] }} \geq \tilde C $$

and replacing $\varepsilon$ by $\delta^{1/\gamma} \varepsilon '$ with $\varepsilon '$ sufficiently small, we should require that $me^{-2}$ should be large enough.

\end{rema}

\begin{rema}

Nash's inequality which was used in by Toscani  \cite{T} says that for all $h$:

$$[ \int | h(v)|^2 dv ]^{1+2/3} \leq C \| h\|^{4/3}_{L^1} \| \nabla h\|^2_{L^2}.$$

In our case, $h =<v>^{\gamma /2} <v>^\alpha f$, so we see that we need a moment estimate on $f$. That estimate, see Villani \cite{VT}, grows up linearly in time, and so we get a bad estimate.

We can also use the result of Desvillettes-Villani \cite{DV} Lemma 7 on Page 43: it says that for any $h$ smooth, for all $\beta >0$, for all $\delta >0$, we have
$$\int h^2 <v>^{2\beta} \leq \delta \int | \nabla h |^2 +C_\delta [ \int h<v>^{5\beta /2} ]^2.$$
Again the choice $h= <v>^{\gamma /2} g$ leads to ask for a value of $\gamma$ close to zero
$$\int h^2 <v>^{2\beta}  \leq C (\int | \nabla h |^2 )^{3/5} [ \int h <v>^{5\beta /2} )^{4/5},$$
using it, we can show that, setting
$$\tilde \gamma = - 1-9/4 \gamma,$$
one has
$$\int | \nabla <v>^{\gamma /2} g |^2    \geq C \| g\|^{10/3}_{L^2} M_{\tilde \gamma} (t) ^{-4/3}.$$
The corresponding additive inequality
$$\int h^2 <v>^{2\beta} \leq C \delta^{5/3} \int |\nabla h|^2 + C \delta^{-5/2} [ \int h <v>^{5\beta /2} ]^2$$
gives
$$\| <v>^{\gamma /2} g \|^2_{\dot H_1}  \geq C\delta^{-5/3} X  -C\delta^{-25/3} [M_{\tilde\gamma}(t) ]^2,$$
and going back to  our  differential inequality (\ref{w6}), we get
$${d\over{dt}} X \leq - \| <v>^{\gamma /2} g \|^2_{\dot H_1}  \bigg\{  C_{coer} - C \varepsilon^{3+\gamma} (1+\varepsilon^2)^\alpha E^{1\over 3}  X^{1\over 3} \bigg\} +[CC_{coer}   +m C \varepsilon^\gamma ] X.$$

Assume that at time $t$, we have
$$ C_{coer} - C \varepsilon^{3+\gamma} (1+\varepsilon^2)^\alpha E^{1\over 3}  X^{1\over 3}  \geq 0,$$
that is
\begin{equation}\label{cond-1}
C \varepsilon^{3+\gamma} (1+\varepsilon^2)^\alpha E^{1\over 3}  X^{1\over 3}  \leq  C_{coer}.
\end{equation}
Then using the above inequality, we have
$${d\over{dt}} X \leq - \bigg\{  C\delta^{-5/3} X  -C\delta^{-25/3} [M_{\tilde\gamma}(t) ]^2 \bigg\}  \bigg\{  C_{coer} - C \varepsilon^{3+\gamma} (1+\varepsilon^2)^\alpha E^{1\over 3}  X^{1\over 3} \bigg\} +[CC_{coer}   +m C \varepsilon^\gamma ] X.$$
Now we note that this is also:
$${d\over{dt}} X \leq -   C\delta^{-5/3} X   \bigg\{  C_{coer} - C \varepsilon^{3+\gamma} (1+\varepsilon^2)^\alpha E^{1\over 3}  X^{1\over 3} \bigg\} +   C\delta^{-25/3} [M_{\tilde\gamma}(t) ]^2 \bigg\{ 2 C_{coer} - C \varepsilon^{3+\gamma} (1+\varepsilon^2)^\alpha E^{1\over 3}  X^{1\over 3} \bigg\}  $$
$$ +[CC_{coer}   +m C \varepsilon^\gamma ] X,$$
then
$${d\over{dt}} X \leq -   C\delta^{-5/3} X   \bigg\{  C_{coer} - C \varepsilon^{3+\gamma} (1+\varepsilon^2)^\alpha E^{1\over 3}  X^{1\over 3} \bigg\} -   C\delta^{-25/3} [M_{\tilde\gamma}(t) ]^2 C \varepsilon^{3+\gamma} (1+\varepsilon^2)^\alpha E^{1\over 3}  X^{1\over 3}  $$
$$+ C\delta^{-25/3} [M_{\tilde\gamma}(t) ]^2  C_{coer} +[CC_{coer}   +m C \varepsilon^\gamma ] X.$$
We see that we still have trouble with the growth rate of $M_{\tilde \gamma} (t)$.

\end{rema}

\vspace{2mm}

\section{The case $\gamma =-3$: weighted energy estimates }
\setcounter{equation}{0}

We adapt the proof given in the previous Section 4, by taking below $\gamma =-3$.

Again setting $g=<v>^{\alpha} f$, we have

$${d\over{dt}} \| g\|^2_{L^2} =
  -  \int_v \int_{v_*} a_{ij} f_* \partial_{v_j}  g \partial_{v_i} g  + \int_v \int_{v_*} a_{ij}  \partial_{v_{*j}} f_*  g \partial_{v_i} g $$
\begin{equation}\label{w1}
+ \int_v\int_{v_*} a_{ij} f_* f   \partial_{v_i} g ~\partial_{v_j}  <v>^{\alpha}
\end{equation}
$$- \int_v \int_{v_*} a_{ij} f_* g \partial_{v_j} f \partial_{v_i} <v>^{\alpha}  + \int_v \int_{v_*} a_{ij}  fg \partial_{v_{*j}} f_*  )  \partial_{v_i} <v>^{\alpha}$$

$$= -I + II    + III  -IV +V.$$

We still have 
\begin{equation}
I \geq \frac{C_{coer}}{2}  \int_v  | \nabla_v [ <v>^{\gamma /2} g ] |^2 -CC_{coer} \int_v  g^2.
\end{equation}

For the term $II$, we have

$$II =  \int_v \int_{v_*} a_{ij}  \partial_{v_{*j}} f_*  g \partial_{v_i} g  =  {1\over 2} \int_v \int_{v_*} a_{ij}  \partial_{v_{*j}} f_*   \partial_{v_i} g^2  = 4\pi \int_v <v>^{2\alpha} f^3 dv$$

and therefore, similarly as in Section 4, we get, with $\alpha \geq -1 - 3/2 \gamma$

\begin{equation}
II \leq   [ \int <v>^2 f ]^{1\over 3}  \| <v>^\alpha f\|^{2\over 3}_{L^2}  \int | \nabla (<v>^{\alpha + \gamma /2} f )|^2 .
\end{equation}

Next, recalling the term $III$ which is given by

$$III = \int_v\int_{v_*} a_{ij} f_* f~   \partial_{v_i} g ~\partial_{v_j}  <v>^{\alpha} ,$$

we have immediately that

$$  III  = \int_v\int_{v_*} a_{ij} f_* f^2~ \partial_{v_i}  <v>^{\alpha} \partial_{v_j}  <v>^{\alpha} + \int_v\int_{v_*} a_{ij} f_* f \partial_{v_i} f <v>^{\alpha} \partial_{v_j}  <v>^{\alpha}
$$
$$=   \int_v\int_{v_*} a_{ij} f_* f^2~ \partial_{v_i}  <v>^{\alpha} \partial_{v_j}  <v>^{\alpha} +
\frac{1}{2} \int_v\int_{v_*} a_{ij} f_*  \partial_{v_i} f^2 <v>^{\alpha} \partial_{v_j}  <v>^{\alpha}
$$
$$=  \frac{1}{2} \int_v\int_{v_*} a_{ij} f_* f^2~ \partial_{v_i}  <v>^{\alpha} \partial_{v_j}  <v>^{\alpha} -
\frac{1}{2} \int_v\int_{v_*} \partial_{v_i}a_{ij} f_* f^2 <v>^{\alpha} \partial_{v_j}  <v>^{\alpha}
$$
$$ -\frac{1}{2} \int_v\int_{v_*}a_{ij} f_* f^2 <v>^{\alpha} \partial_{v_i} \partial_{v_j}  <v>^{\alpha},$$

Since 
$$\partial_{v_i} <v>^\alpha = \partial_{v_i} (1+|v|^2 )^{\alpha /2} = {\alpha \over 2} (1+|v|^2 )^{(\alpha -2) /2 } 2v_i \mbox{ and }| \partial_{v_i} <v>^\alpha | \leq C (1+|v|^2 )^{(\alpha -1) /2 } \leq C <v>^{\alpha -1},$$

we obtain

$$III \leq  C \int_v\int_{v_*} |v-v_*|^{\gamma +2} f_*f^2  <v>^{2\alpha -2} + \int_v\int_{v_*} |v-v_*|^{\gamma +1} f_*f^2 <v>^{2\alpha -1} +\int_v\int_{v_*} |v-v_*|^{\gamma +2} f_* f^2 <v>^{2\alpha -2}$$

Then, we have for any $\varepsilon >0$ fixed

$$ III \leq III_\varepsilon + \overline{III}_\varepsilon$$

where

$$\overline{III}_\varepsilon = C \int_v\int_{v_*}1_{|v-v_*|\geq \varepsilon} |v-v_*|^{\gamma +2} f_*f^2  <v>^{2\alpha -2} + \int_v\int_{v_*} 1_{|v-v_*|\geq \varepsilon} |v-v_*|^{\gamma +1} f_*f^2 <v>^{2\alpha -1} $$
$$+\int_v\int_{v_*}1_{|v-v_*|\geq \varepsilon}  |v-v_*|^{\gamma +2} f_* f^2 <v>^{2\alpha -2}$$

and

$$III_\varepsilon = C \int_v\int_{v_*}1_{|v-v_*|\leq \varepsilon} |v-v_*|^{\gamma +2} f_*f^2  <v>^{2\alpha -2} + \int_v\int_{v_*} 1_{|v-v_*|\leq \varepsilon} |v-v_*|^{\gamma +1} f_*f^2 <v>^{2\alpha -1} $$
$$+\int_v\int_{v_*}1_{|v-v_*|\leq \varepsilon}  |v-v_*|^{\gamma +2} f_* f^2 <v>^{2\alpha -2} .$$

For $\overline{III}_\varepsilon $, one has

$$\overline{III}_\varepsilon \leq  C \int_v m\varepsilon^{\gamma +2} f^2  <v>^{2\alpha -2} + \int_v m\varepsilon^{\gamma +1} f^2 <v>^{2\alpha -1}+\int_vm \varepsilon^{\gamma +2} f^2 <v>^{2\alpha -2}$$

and thus

$$\overline{III}_\varepsilon \leq  C  m [ \varepsilon^{\gamma +2}  + \varepsilon^{\gamma +1} ] \| g\|^2_{L^2} .$$

For $III_\varepsilon$, splitting over the sets $f_*\leq f$ and $f\leq f_*$, we obtain that

$$III_\varepsilon \leq A+B$$

where

$$A \leq C \int_v\int_{v_*}1_{|v-v_*|\leq \varepsilon} |v-v_*|^{\gamma +2} f^3  <v>^{2\alpha -2} + \int_v\int_{v_*} 1_{|v-v_*|\leq \varepsilon} |v-v_*|^{\gamma +1} f^3 <v>^{2\alpha -1} $$
$$+\int_v\int_{v_*}1_{|v-v_*|\leq \varepsilon}  |v-v_*|^{\gamma +2} f^3 <v>^{2\alpha -2}$$

and

$$B \leq C \int_v\int_{v_*}1_{|v-v_*|\leq \varepsilon} |v-v_*|^{\gamma +2} f_*^3  <v>^{2\alpha -2} + \int_v\int_{v_*} 1_{|v-v_*|\leq \varepsilon} |v-v_*|^{\gamma +1} f_*^3 <v>^{2\alpha -1} $$
$$+\int_v\int_{v_*}1_{|v-v_*|\leq \varepsilon}  |v-v_*|^{\gamma +2} f_*^3 <v>^{2\alpha -2} .$$

We have

$$A \leq C  \varepsilon^{3+2+\gamma} \int_v f^3  <v>^{2\alpha -2} + \varepsilon^{3+\gamma +1} \int_v   f^3 <v>^{2\alpha -1} + \varepsilon^{3+\gamma +2} \int_v  f^3 <v>^{2\alpha -2}$$
and we see immediately that

$$A \leq C [  \varepsilon^{3+2+\gamma} + \varepsilon^{3+\gamma +1}  + \varepsilon^{3+\gamma +2} ] II$$
and therefore

$$A \leq C [  \varepsilon^{3+2+\gamma} + \varepsilon^{3+\gamma +1}  + \varepsilon^{3+\gamma +2} ]  [ \int <v>^2 f ]^{1\over 3}  \| <v>^\alpha f\|^{2\over 3}_{L^2}  \int | \nabla (<v>^{\alpha + \gamma /2} f )|^2 .$$

For $B$, the same arguments leads to

$$B  \leq C(1+\varepsilon^2)^\alpha [  \varepsilon^{3+2+\gamma} + \varepsilon^{3+\gamma +1}  + \varepsilon^{3+\gamma +2} ]  [ \int <v>^2 f ]^{1\over 3}  \| <v>^\alpha f\|^{2\over 3}_{L^2}  \int | \nabla (<v>^{\alpha + \gamma /2} f )|^2 $$

and thus

$$III_\varepsilon  \leq C(1+\varepsilon^2)^\alpha [  \varepsilon^{3+2+\gamma} + \varepsilon^{3+\gamma +1}  + \varepsilon^{3+\gamma +2} ]  [ \int <v>^2 f ]^{1\over 3}  \| <v>^\alpha f\|^{2\over 3}_{L^2}  \int | \nabla (<v>^{\alpha + \gamma /2} f )|^2 .$$

In conclusion, we get

$$III  \leq C(1+\varepsilon^2)^\alpha [  \varepsilon^{3+2+\gamma} + \varepsilon^{3+\gamma +1}  + \varepsilon^{3+\gamma +2} ]  [ \int <v>^2 f ]^{1\over 3}  \| <v>^\alpha f\|^{2\over 3}_{L^2}  \int | \nabla (<v>^{\alpha + \gamma /2} f )|^2 $$
$$ + C  m [ \varepsilon^{\gamma +2}  + \varepsilon^{\gamma +1} ] \| g\|^2_{L^2} .$$

The same arguments can be applied to all other terms, and thus we get

$$II +|III| +| IIV|+ |V|   \leq C \bigg\{ 1+ (1+\varepsilon^2)^\alpha [  \varepsilon^{3+2+\gamma} + \varepsilon^{3+\gamma +1}  + \varepsilon^{3+\gamma +2} ] \bigg\} [ \int <v>^2 f ]^{1\over 3}  \| <v>^\alpha f\|^{2\over 3}_{L^2}  \int | \nabla (<v>^{\alpha + \gamma /2} f )|^2 $$
$$ + C  m [ \varepsilon^{\gamma +2}  + \varepsilon^{\gamma +1} ] \| g\|^2_{L^2} .$$

Note the difference when $\gamma >-3$: in that case, the constant is small in the first term, while here for $\gamma =-3$, we have a constant which is close to $1$, for small $\varepsilon$.

We let $O(\varepsilon )$ for the first function and $\tilde O (1/\varepsilon )$ for the second one to get

$$II +|III| +| IIV|+ |V|   \leq C \bigg\{ 1+ O(\varepsilon ) \bigg\} [ \int <v>^2 f ]^{1\over 3}  \| <v>^\alpha f\|^{2\over 3}_{L^2}  \int | \nabla (<v>^{\alpha + \gamma /2} f )|^2 $$
$$ + C  m \tilde O({1\over\varepsilon}) \| g\|^2_{L^2} .$$

By combining the above estimations, the final conclusion is that

\begin{equation}
 {d\over{dt}} \| g\|^2_{L^2} +  C_{coer}  \int_v  | \nabla_v [ <v>^{\gamma /2} g ] |^2 \leq
\end{equation}
$$CC_{coer} \| g\|^2_{L^2}   + C \bigg\{ 1+ O(\varepsilon ) \bigg\} [ \int <v>^2 f ]^{1\over 3}  \| <v>^\alpha f\|^{2\over 3}_{L^2}  \int | \nabla (<v>^{\alpha + \gamma /2} f )|^2 $$
$$ + C  m \tilde O({1\over\varepsilon}) \| g\|^2_{L^2}. $$

Again setting $E= \int <v>^2 f dv = m+e$, which is bounded uniformly in time, let
$$X = \| g\|^2_{L^2},   \hspace{5mm}  \| <v>^{\gamma /2} g \|^2_{\dot H_1} = \int | \nabla (<v>^{\gamma /2} g )|^2, $$
then we have obtained
 \begin{equation}
 {d\over{dt}} X \leq - \| <v>^{\gamma /2} g \|^2_{\dot H_1}  \bigg\{  C_{coer} - C (1+O(\varepsilon ))  E^{1\over 3}  X^{1\over 3} \bigg\} +[CC_{coer}   +m C \tilde O ({1\over\varepsilon}) ] X .
\end{equation}

At this point, we can use the arguments of Section 4. Assume that at time $t$, we have for some $\delta \geq 0$
$$ C_{coer} - C (1+O(\varepsilon )) E^{1\over 3}  X^{1\over 3}  \geq \delta, $$
that is
\begin{equation}
C (1+O(\varepsilon )) E^{1\over 3}  X^{1\over 3}  \leq  C_{coer} -\delta,
\end{equation}
which is to be used for $\delta < C_{coer}$.

Then we obtain, using again Pitt's inequality (\ref{w8}) and (\ref{w10}),
\begin{equation}\label{radja-1}
{d\over{dt}} X \leq - C_{pitt} C m^{11/2}  e^{-7/2} \delta +[CC_{coer}   +m \tilde O ({1\over\varepsilon}) ] X =F(X).
\end{equation}

Recall that we have assumed 
$$\delta <  C_{coer}$$
and that we want
$$C (1+O(\varepsilon )) E^{1\over 3}  X^{1\over 3}  \leq  C_{coer} -\delta ,$$
that is
$$
  X\leq  [ C_{coer} -\delta  ]^3 C^{-3} (1+O(\varepsilon ))^{-3} / E \equiv \tilde X .
$$
Now, let $X_{eq}$ be the zero of the function $F$ defined in (\ref{radja-1}).
Then assume that
\begin{equation}
X(0) \leq \bar X\equiv \min \lbrace \tilde X, X_{eq} \rbrace,
\end{equation}
then in view of the form of the differential inequality and the behaviour of function $F$, it follows that for all $t>0$:
 $$X(t) \leq \bar X.$$
Thus we have obtain a global bound for the weighted $L^2$ norm of $f$, uniformly in time, that is, we get the second part of Theorem \ref{Th}, for the specific case $\gamma =-3$.

\section{The case $-3< \gamma <-2$: local estimates}
\setcounter{equation}{0}

The energy estimate in Section 3 holds for $\gamma \geq -2$. For the case $\gamma \in (-3,-2)$, we recall (\ref{e1})
$${d\over{dt}} {1\over 2} \| f\|^2_{L^2} +\int\int a_{ij} (v-v_*) f_* \partial_{v_i} f\partial_{v_j} f =
{1\over 2}\int\int a_{ij} (v-v_*) \partial_{v_{*j}} f_* \partial_{v_i} f^2,
$$
and from (\ref{e2}) and (\ref{e7}) we have
$$
\int\int a_{ij} (v-v_*) f_* \partial_{v_i} f\partial_{v_j} f
\geq \frac{C_{coer}}{2}  \int  |\nabla_v (<v>^{\gamma/2}f)|^2 dv
- C \int f^2 dv.
$$
Next, we estimate
$${1\over 2}\int\int a_{ij} (v-v_*) \partial_{v_{*j}} f_* \partial_{v_i} f^2 = (\gamma +3) \int_{v_*} \int_v |v-v_*|^\gamma f_* f^2.  $$
The problem is that it looks like a $L^3$ norm, but at that point we need a $L^1$ weighted estimation. Up to now these bounds grow linearly in time \cite{VT} and so are not enough.

Let us fix a positive function of time $\phi (t)$.
We split $A$ into two terms (forgetting the positive constant $\gamma +3$ in front of $A$)
$$A =A_1 +A_2,$$
where
$$A_1 =  \int_{v_*} \int_v 1_{|v-v_*| \geq \phi (t)} |v-v_*|^\gamma f_* f^2 $$
and
$$A_2 =  \int_{v_*} \int_v 1_{|v-v_*| \leq \phi (t)} |v-v_*|^\gamma f_* f^2. $$

For $A_1$, since $\gamma <0$, we have
\begin{equation}\label{l1}
A_1 \leq  \int_{v_*} \int_v \phi (t)^\gamma f_* f^2 \leq m \phi (t)^\gamma \| f\|^2_{L^2}.
\end{equation}

For $A_2$, we split again according to whether or not $f\leq f_*$ to get
\begin{equation}\label{l2}
A_2  \leq 2 \int_{v_*} \int_v 1_{|v-v_*| \leq \phi (t)} |v-v_*|^\gamma f^3 \leq 2 C \phi (t)^{3+\gamma} \| f\|^3_{L^3}.
\end{equation}

Next, we are going to work on $\| f\|_{L^3}$: if $f= f_1f_2f_3$, then, with ${1\over 3} = {1\over p_1}+{1\over p_2} +{1\over p_3}$, we have
$$
\| f\|_{L_3} \leq \| f_1\|_{L_{p_1}} \| f_2\|_{L_{p_2}} \| f_2\|_{L_{p_2}}.
$$
Now write
$$1 = {1\over p_1} +{2\over p_2} + {6\over p_3},$$

$$f =f^{1\over p_1} \cdot f^{2\over p_2} \cdot f^{6\over p_3}$$
$$=f^{1\over p_1} <v>^{-3 \gamma \over p_3}\cdot f^{2\over p_2} \cdot (<v>^{\gamma \over 2} f)^{6\over p_3}$$
$$=( <v>^{-3 p_1\gamma \over p_3} f) ^{1\over p_1} \cdot f^{2\over p_2} \cdot (<v>^{\gamma\over 2} f)^{6\over p_3}$$
$$:=f_1 \cdot f_2\cdot f_3$$
with evident notations.
We make the choice $p_1 =p_2 =p_3 =9$ for reasons linked to Sobolev inequality. Then we get
$$\|f \|_{L^3}^3 \lesssim  \| f_1\|^3_{L_{p_1}} \| f_2\|^3_{L_{p_2}} \| f_2\|^3_{L_{p_2}} $$
$$\lesssim  [ \int <v>^{-3\gamma} f ]^{1\over 3}  \cdot [  \int f^2 ]^{1\over 3} \cdot [ \int (<v>^{\gamma\over 2} f )^6 ]^{1\over 3}.$$
Sobolev inequality tells us that
$$[ \int (<v>^{\gamma\over 2} f )^6 ]^{1\over 3} \lesssim \int | \nabla (<v>^{\gamma\over 2} f )|^2 .$$
Finally, we have obtained
$$\|f \|_{L^3}^3 \leq C  [ \int <v>^{-3\gamma} f ]^{1\over 3}  \| f\|^{2\over 3}_{L^2}  \int | \nabla (<v>^{\gamma\over 2} f )|^2 ,$$
and thus (\ref{l2}) becomes
$$A_2 \leq C \phi (t)^{3+\gamma}  [ \int <v>^{-3\gamma} f ]^{1\over 3}  \| f\|^{2\over 3}_{L^2}  \int | \nabla (<v>^{\gamma\over 2} f )|^2 .$$

From Villani \cite{VT} (Appendix B), in our case $\gamma\geq -3$, i.e., $-3\gamma \leq 9$,  we have
\begin{equation}\label{moment}
M_{-3\gamma} := \int <v>^{-3\gamma} f \leq C (1+t),
\end{equation}
then finally
$$A_2 \leq C \phi (t)^{3+\gamma}  (1+t)^{1\over 3}  \| f\|^{2\over 3}_{L^2}  \int | \nabla (<v>^{\gamma\over 2} f )|^2 .$$

Now we choose $\phi (t)$ such that (for some $\varepsilon$ fixed)
$$\phi (t)^{3+\gamma } (1+t)^{1\over 3} \leq \varepsilon,$$
that is
$$ \phi (t) \leq \varepsilon^{1\over{3+\gamma }} (1+t)^{-1\over {3 (3+\gamma )}} ,$$
or, for simplicity, we just choose
\begin{equation}\label{l3}
\phi (t) =\varepsilon^{1\over{3+\gamma }} (1+t)^{-1\over {3 (3+\gamma )}} .
\end{equation}
With this choice, we get that
$$A_2 \leq C \varepsilon \| f\|^{2\over 3}_{L^2}  \int | \nabla (<v>^{\gamma\over 2} f )|^2 .$$
Recall (\ref{l1}) to get also
$$
A_1 \leq m \phi (t)^\gamma \| f\|^2_{L^2} = m \varepsilon^{\gamma \over{3+\gamma }} (1+t)^{-\gamma \over {3 (3+\gamma )}}  \| f\|^2_{L^2} .
 $$

In conclusion, we get:
\begin{equation}\label{l4}
A \leq m \varepsilon^{\gamma \over{3+\gamma }} (1+t)^{-\gamma \over {3 (3+\gamma )}}  \| f\|^2_{L^2}   +   C \varepsilon \| f\|^{2\over 3}_{L^2}  \int | \nabla (<v>^{\gamma\over 2} f )|^2 .
\end{equation}

Fix $t$, and optimize w.r.t. $\varepsilon$. The above term is of the form
$$\varepsilon^{\gamma \over{3+\gamma }} B + \varepsilon D$$
with
$$B = m(1+t)^{-\gamma \over {3 (3+\gamma )}}  \| f\|^2_{L^2}, \hspace{5mm}
D =  C  \| f\|^{2\over 3}_{L^2}  \int | \nabla (<v>^{\gamma\over 2} f )|^2. $$
We have equality in (\ref{l4}) if
$$\varepsilon^{\gamma \over{3+\gamma }} B = \varepsilon D,$$
that is
$$\varepsilon = B^{ {3+\gamma }\over{3}  } D^{ - {{3+\gamma }\over{3}} }.$$
With this value, we get from  (\ref{l4}) that
$$A \leq 2 \varepsilon D = 2 B^{ {3+\gamma }\over{3}  } D^{ - {{3+\gamma }\over{3}} } D= 2 B^{ {3+\gamma }\over{3}  } D^{ - {{\gamma }\over{3}} },$$
 that is
$$A \leq 2  [  m(1+t)^{-\gamma \over {3 (3+\gamma )}}  \| f\|^2_{L^2}  ]^{ {3+\gamma }\over{3}  }  \ [  C  \| f\|^{2\over 3}_{L^2}  \int | \nabla (<v>^{\gamma\over 2} f )|^2  ]^{ - {{\gamma }\over{3}} }$$

$$\leq C  (1+t)^{ - \gamma\over 9}  [   \int | \nabla (<v>^{\gamma\over 2} f )|^2  ]^{ - {{\gamma }\over{3}} } \| f\|_{L^2}^{ {18+4\gamma }\over{9}}.$$

  We choose another $\varepsilon$ and write

$$A \leq C  (1+t)^{ - \gamma\over 9}  [   \int | \nabla (<v>^{\gamma\over 2} f )|^2  ]^{ - {{\gamma }\over{3}} } \| f\|_{L^2}^{ {18+4\gamma }\over{9}}$$

$$\leq C  \varepsilon^{\frac{\gamma}{3}} (1+t)^{ - \gamma\over 9}  \varepsilon^{-\frac{\gamma}{3}} [   \int | \nabla (<v>^{\gamma\over 2} f )|^2  ]^{ - {{\gamma }\over{3}} } \| f\|_{L^2}^{ {18+4\gamma }\over{9}},$$
use Young's inequality for product of the first two factors (ie without the $L^2$ norm) with $p = -3/\gamma$ and $p' = 3/(3+\gamma )$:
$$
A \leq \bigg(  \frac{\varepsilon}{p}  \int | \nabla (<v>^{\gamma\over 2} f )|^2    + C \varepsilon^{\frac{\gamma}{3+\gamma}} (1+t)^{-\frac{\gamma}{3(3+\gamma)}}  \bigg) \ \| f\|_{L^2}^{ {18+4\gamma }\over{9}}.
$$

Combine the above estimates, we get a differential inequality
\begin{equation}\label{l5}
{d\over{dt}}  \| f\|^2_{L^2}   \leq
-\bigg(c_{coer} - \frac{\varepsilon}{p} \| f\|_{L^2}^{ {18+4\gamma }\over{9}}\bigg)
\int | \nabla (<v>^{\gamma\over 2} f )|^2
 + C \| f\|_{L^2} + C (1+t)^{-\frac{\gamma}{3(3+\gamma)}}
 \| f\|_{L^2}^{ {18+4\gamma }\over{9}}.
\end{equation}

We see that we have trouble with the growth rate of coefficient, due to lack of uniform in time bound of the moment in (\ref{moment}). The differential inequality (\ref{l5}) yields a local estimate of a weak solution in this case and thus we have the following weaker conclusion

\begin{prop}
Let $\gamma\in (-3,-2)$. Let the initial data $f_0\in L^2(\RR^3)$, then we have a local in time a priori estimate \eqref{l5} on a weak solution in $L^2$.
\end{prop}

\begin{rema}
 We emphasize here that when $\gamma\in (-3,-2)$, the a priori estimate in $L^2$ is only local, unless we can get a better moment estimates in $L^1$,that is, uniformly bounded w.r.t. time.
\end{rema}

\section{Proof of Proposition \ref{prop-radja}}
\setcounter{equation}{0}

Multiplying the Landau equation by $f$ and integrating, we have

$${d\over{dt}} {1\over 2} \| f\|^2_{L^2} +\int\int a_{ij} (v-v_*) f_* \partial_{v_i} f\partial_{v_j} f =
{1\over 2}\int\int a_{ij} (v-v_*) \partial_{v_{*j}} f_* \partial_{v_i} f^2.
$$
Moreover, as usual now, we have also (note here that we keep the weight on the second term on the r.h.s.)
$$
\int\int a_{ij} (v-v_*) f_* \partial_{v_i} f\partial_{v_j} f
\geq \frac{C_{coer}}{2}  \int  |\nabla_v (<v>^{\gamma/2}f)|^2 dv
- C \int <v>^{\gamma -2} f^2 dv.
$$
and
$${1\over 2}\int\int a_{ij} (v-v_*) \partial_{v_{*j}} f_* \partial_{v_i} f^2 = (\gamma +3) \int_{v_*} \int_v |v-v_*|^\gamma f_* f^2,  $$

All in all, we have

$${d\over{dt}} {1\over 2} \| f\|^2_{L^2} + \frac{C_{coer}}{2}  \int  |\nabla_v (<v>^{\gamma/2}f)|^2 dv \leq (\gamma +3) \int_{v_*} \int_v |v-v_*|^\gamma f_* f^2 + C \int <v>^{\gamma -2} f^2 dv.$$

Define the first term on the r.h.s. as $NLT$ (non linear term), that is

$$NLT = (\gamma +3) \int_{v_*} \int_v |v-v_*|^\gamma f_* f^2.$$
From now on, we will omit or abbreviate any non important constant. For any $\varepsilon >0$, we can write

$$NLT \leq C  \int_v f^2(v) \bigg\{ \int_{v_*} |v-v_*|^\gamma 1_{|v-v_*|\leq \varepsilon} f_* + \int_{v_*} |v-v_*|^\gamma 1_{|v-v_*|\geq \varepsilon} f_* \bigg\} .$$

Then, we can use classical estimations on the truncated Riez potentials, see \cite{ZIE} for example, involving the usual maximal function $Mf(v)$ to get

$$NLT \leq C  \int_v f^2(v) \bigg\{ \varepsilon^{3+\gamma} Mf(v) + \varepsilon^\gamma m \bigg\}.$$

Fixing $v$, we optimize w.r.t. $\varepsilon$ to find that

$$NLT \leq C m^{1+\gamma /3} \int_v f^2(v) Mf(v)^{-\gamma /3} \leq C C m^{1+\gamma /3}  \int Mf (v)^{2-\gamma /3} \leq C m^{1+\gamma /3} \int f^{2-\gamma /3},$$
by using our assumption on the values of $\gamma$. 

Let $q$ be defined by $q= 2- \gamma /3 = {{6-\gamma}\over 3} >1$. Note that we have also $q_1 = {q\over 2} = {{6-\gamma}\over 6} >1$. The conjugate exponent is given by $q_1' = {{6-\gamma}\over{-\gamma}}$. We can then use Holder inequality together with the fact that $\gamma -2 \leq -3$ to get

$$C \int <v>^{\gamma -2} f^2 dv \leq C \bigg\{ \int f^{{6-\gamma}\over 3} \bigg\}^{6\over{6-\gamma}} .$$

Using the conservation of mass, again skipping all constants, we have obtained

\begin{equation}\label{radja-new-1}{d\over{dt}} {1\over 2} \| f\|^2_{L^2} + \frac{C_{coer}}{2}  \int  |\nabla_v (<v>^{\gamma/2}f)|^2 dv \leq C_1 +C_2 \int f^{{6-\gamma}\over 3} dv .
\end{equation}

Now the idea is this: we want to control the l.h.s. term by the r.h.s, and so we will use some Nash Gagliardo Nirenberg type inequalities, \cite {TAY} for example.

We have a slight issue connected to moments (because on the l.h.s., we have only some control of a negative power weight in Sobolev space), but let's forget this point for the moment. Firstly recall that (we are using homogeneous Sobolev spaces) $\dot H_s \subset_c L^m$ for $0<s<3/2$ and $m = {6\over{3-2s}}$. We want to choose $m=q = {{6-\gamma}\over 3}$. This gives the value of $s$ as $s= {{3\gamma}\over{2\gamma -12}}$. Note that we have $0<s<1$. 

If this is the case, it follows that $\| f\|^q_{L^q} \leq C \| f\|^q_{\dot H_s}$. 

On the other hand, by using classical ideas for proving Nash inequality (Fourier transform, optimizing for small and big frequencies), one can show that (for $s<1$ which is the case here)

$$\| f\|_{\dot H_s} \leq C m^{1- {1\over 5} (3+2s)} \bigg\{ \int | \nabla f\|^2  \bigg\}^{{1\over{10}} (3+2s)}$$

and thus

$$\| f\|_{L^q}^q \leq C m^{[1- {1\over 5} (3+2s)]q} \bigg\{ \int | \nabla f\|^2  \bigg\}^{{1\over{10}} (3+2s)q} .$$

A little computation shows that $\mu \equiv {1\over{10}} (3+2s)q = {1\over 5} [3-\gamma ]$ which gives $\mu <1$ iff $\gamma > -2$. Then (up to the control of weights), we can absorb the r.h.s by the l.h.s in inequality \eqref{radja-new-1}.

Now to get everything rigorous, and in particular to take care of the loss of weights appearing on the l.h.s, we need to interpolate with a weighted $L^1$ space the r.h.s. of \eqref{radja-new-1} (as well we can also use some improved type Nash inequalities).

Starting with a fixed $\varepsilon >0$, we look for $\alpha \in (0,1)$ such that $q= \alpha .1 +(1-\alpha ) (q+\varepsilon )$. We find that
$$\alpha = {\varepsilon\over{q+\varepsilon -1}} \mbox{ and } 1-\alpha ={{q-1}\over{q+\varepsilon -1}}.$$

It follows that
$$\int f^q = \int f^{\alpha .1 +(1-\alpha ) (q+\varepsilon )}= \int f^{\alpha} f^{ (1-\alpha)(q+\varepsilon )} = \int  [<v>^{ {{-\gamma}\over{2\alpha}} (1-\alpha)(q+\varepsilon )  }f]^\alpha [<v>^{\gamma /2}f]^{ (1-\alpha)(q+\varepsilon )} $$
and using Holder inequality, we get

$$\int f^q \leq \bigg\{ M_{ {{-\gamma}\over{2\alpha}} (1-\alpha)(q+\varepsilon )  } \bigg\}^{\alpha} . \bigg\{ \int [ <v>^{\gamma /2} f ]^{q+\varepsilon} \bigg\}^{1-\alpha},$$
which upon using another small $\tilde \varepsilon >0$, yields

$$\int f^q \leq C(\tilde \varepsilon ) M_{ {{-\gamma}\over{2\alpha}} (1-\alpha)(q+\varepsilon )  }  (t)+ \tilde \varepsilon \int [ <v>^{\gamma /2} f ]^{q+\varepsilon} .$$

Set $q_\varepsilon =q+\varepsilon$, $\gamma_\varepsilon = \gamma -3\varepsilon$ and $s_\varepsilon = {{2\gamma_\varepsilon}\over{2\gamma_\varepsilon -12}}$.

Then, we still have

$$\| g\|_{L^{q_\varepsilon}}^{q_\varepsilon} \leq C m^{[1- {1\over 5} (3+2s_\varepsilon)]q_\varepsilon} \bigg\{ \int | \nabla g\|^2  \bigg\}^{{1\over{10}} (3+2s_\varepsilon)q_\varepsilon} .$$

Then note that $\mu_\varepsilon = {1\over{10}} (3+2s_\varepsilon)q_\varepsilon = {1\over 5} [ 3-\gamma_\varepsilon ]= {1\over 5} [ 3-\gamma +3\varepsilon ]$. Since we have assumed $\gamma >-2$, one obtains that $\mu_\varepsilon \leq 1$ when choosing any $\varepsilon$ such that $0< \varepsilon \leq {{2+\gamma}\over 3}$. Therefore we choose exactly $\varepsilon = {{2+\gamma}\over 3}$. 

With this value of $\varepsilon$, we find that $\alpha = {{2+\gamma}\over 5}$, $1-\alpha ={{3-\gamma}\over 5}$, $q+\varepsilon = {8\over 3}$, and thus all in all, we find that

$$\int f^q \leq C (\tilde\varepsilon ) M_{ {-4\gamma (3-\gamma )}\over{3(2+\gamma )}} (t) + C(\tilde \varepsilon ,m) + \tilde \varepsilon C(m) \int |\nabla [<v>^{\gamma /2} f ]|^2 ,$$

and therefore choosing $\tilde \varepsilon $ small enough, we find that, going back to our estimation inequality
 
$$ \| f \|^2_{L^2} (t) \leq C  (1+t)^2 $$
which ends the proof of Proposition \ref{prop-radja}.

\begin{rema}

\begin{enumerate}

\item Of course, we get also a Sobolev estimation as well. Moreover, it might be also possible to have direct estimation of the nonlinear term by using Holder inequality, together with the standard Nash's inequality.

\item Note that the growth of this $L^2$ estimate is linked with the moment estimate. One can also get weighted $L^2$ estimate and more generally $L^p$ estimates. For example, one can show that the nonlinear term is estimated by $\int f^{p-\gamma /3}$. But it does not seem to be possible to improve the range of values of $\gamma$. However, working with large $p$ seem to require less moments on $f$.

\item By interpolating also with $L^2$, as in previous sections, one can get also local estimates for all $\gamma >-3$.

\end{enumerate}

\end{rema}

\smallskip

{\bf Acknowledgements:} This work was supported by the Fundamental Research Funds for the Central Universities and National Natural Science Foundation of China (Nos.11171211, 11171212), together with a starting grant from Shanghai Jiao Tong University.

\end{document}